\documentclass[12pt]{article}
\usepackage{amsmath,amssymb,amsthm}
\newcommand{\tmmathbf}[1]{\ensuremath{\boldsymbol{#1}}}

\newtheorem{theorem}{Theorem}[section]
\newtheorem{corollary}[theorem]{Corollary}
\newtheorem{lemma}[theorem]{Lemma}

\newtheorem{definition}[theorem]{Definition}

\numberwithin{equation}{section}

\def\R{\mathcal{R}}

\begin{document}

\title{Conditional expanding bounds for two-variable functions over finite valuation rings}

\author{
Le Quang Ham\thanks{University of Science, Vietnam National University Hanoi
    Email: {hamlaoshi@gmail.com
}}
\and 
    Pham Van Thang\thanks{EPFL, Lausanne, Switzerland. Research partially supported by Swiss National Science Foundation Grants 200020-144531, 200021-137574 and 200020-162884
    Email: {\tt thang.pham@epfl.ch}}
  \and
    Le Anh Vinh\thanks{University of Education, Vietnam National University Hanoi. Research was supported by Vietnam National Foundation for Science and Technology Development grant 101.99-2013.21.
    Email: {\tt vinhla@vnu.edu.vn
}}}
\date{}
\maketitle
\begin{abstract}
In this paper, we use methods from spectral graph theory to obtain some results on the sum-product problem over finite valuation rings $\mathcal{R}$ of order $q^r$ which generalize recent results given by Hegyv\'{a}ri and Hennecart (2013). More precisely, we prove that, for related pairs of two-variable functions $f(x,y)$ and $g(x,y)$, if $A$ and $B$ are two sets in $\mathcal{R}^*$ with $|A|=|B|=q^\alpha$, then 
\[\max\left\lbrace |f(A, B)|, |g(A, B)| \right\rbrace\gg |A|^{1+\Delta(\alpha)},\]
for some $\Delta(\alpha)>0$.
\end{abstract}

\section{Introduction}

Let $\mathbb{F}_q$ be a finite field of $q$ elements where $q$ is an odd prime power. Throughout the paper $q$ will be a large prime power. Let $\mathcal{A}$ be a non-empty subset of a finite field $\mathbb{F}_q$. We consider the sum set
\[ \mathcal{A}+\mathcal{A}:=\{a+b : a, b \in \mathcal{A}\}\]
and the product set
\[\mathcal{A}\cdot \mathcal{A}:=\{a\cdot b : a, b \in \mathcal{A}\}.\]

Let $|\mathcal{A}|$ denote the cardinality of $\mathcal{A}$. Bourgain, Katz and Tao (\cite{bourgain-katz-tao}) showed  that when $ 1\ll |\mathcal{A}| \ll q$ then
$\max(|\mathcal{A}+\mathcal{A}|,|\mathcal{A}\cdot\mathcal{A}|) \gg |\mathcal{A}|^{1+\epsilon},$ for some $\epsilon >0$. This improves the trivial bound $\max\{|\mathcal{A}+\mathcal{A}|, |\mathcal{A}\cdot \mathcal{A}|\}\gg |\mathcal{A}|$. (Here, and throughout, $X \asymp Y$ means that there exist positive constants $C_1$ and $C_2$ such that $C_1Y< X <C_2Y$, and $X \ll Y$ means that there exists $C>0$ such that $X\le  CY$). The precise statement of their result is as follows.

\begin{theorem}[\textbf{Bourgain, Katz and Tao}, \cite{bourgain-katz-tao}] \label{t1}
Let $\mathcal{A}$ be a subset of $\mathbb{F}_q$ such that
\[ q^\delta < |\mathcal{A}| < q^{1-\delta}\]
for some $\delta > 0$. Then one has a bound of the form
\[ \max \left\lbrace|\mathcal{A}+\mathcal{A}|, |\mathcal{A}\cdot\mathcal{A}|\right\rbrace \gg |\mathcal{A}|^{1+\epsilon}\]
for some $\epsilon = \epsilon(\delta) > 0$.
\end{theorem}
Note that the relationship between $\epsilon$ and $\delta$ in Theorem \ref{t1} is difficult to determine. In \cite{his}, Hart, Iosevich, and Solymosi obtained a bound that gives an explicit dependence of $\epsilon$ on $\delta$. More precisely, if $|A + A| = m$ and $|A \cdot A | = n$, then 
\begin{equation}\label{eq:his}
|A|^3 \leq \frac{ c m^2 n |A| }{ q} + c q^{1/2} mn,
\end{equation}
for some positive constant $c$.  Inequality (\ref{eq:his}) implies a non-trivial sum-product estimate when 
$|A|\gg q^{1/2}$.  Using methods from the spectral graph theory, the third listed author \cite{vinh} improved (\ref{eq:his}) and as a result, obtained a better sum-product estimate.     

\begin{theorem}[\textbf{Vinh}, \cite{vinh}]\label{vinh th}
For any set $A\subseteq \mathbb{F}_q$, if $|A+A| = m$, and $|A \cdot A| = n$, then 
\[
|A|^2 \leq \frac{mn|A|}{q} + q^{1/2} \sqrt{ m n } .
\]
\end{theorem}

\begin{corollary}[\textbf{Vinh}, \cite{vinh}]\label{vinh cor}
For any set $A \subseteq \mathbb{F}_q$, we have\\
If $q^{1/2} \ll |A| \ll q^{2/3}$, then 
\[
\max \left\lbrace |A + A | , |A \cdot A | \right\rbrace \gg \frac{ |A|^2 }{q^{1/2}} .
\]
If $|A| \gg q^{2/3}$, then 
\[
\max \left\lbrace |A + A | , |A \cdot A | \right\rbrace \gg ( q |A| )^{1/2} .
\]
\end{corollary}
It follows from Corollary \ref{vinh cor} that if $|A|=p^\alpha$, then
\[\max\left\lbrace |A+A|, |A\cdot A|\right\rbrace\gg |A|^{1+\Delta(\alpha)}, \]
where $\Delta(\alpha)=\min\left\lbrace 1-1/2\alpha, (1/\alpha-1)/2\right\rbrace$.
In the case that $q$ is a prime, Corollary \ref{vinh cor} was proved by Garaev \cite{garaev} using exponential sums.  
 Cilleruelo \cite{cil} also proved related results using dense Sidon sets in finite groups involving $\mathbb{F}_q$ and $\mathbb{F}_q^*:=\mathbb{F}_q\setminus \{0\}$ (see \cite[Section 3]{cil} for more details).

 We note that a variant of Corollary \ref{vinh cor} was considered by Vu \cite{vu}, and the statement is as follows.
\begin{theorem}[\textbf{Vu}, \cite{vu}]\label{vuvu}
Let $P$ be a \textit{ non-degenerate} polynomial of degree $k$ in $\mathbb{F}_q[x, y]$. Then for any $A\subseteq \mathbb{F}_q$, we have 
\[\max\left\lbrace|A+A|, |P(A)|\right\rbrace\gtrsim \min\left\lbrace|A|^{2/3}q^{1/3}, |A|^{3/2}q^{-1/4}\right\rbrace,\]
where we say that a polynomial $P$ is \textit{non-degenerate} if $P$ can not be presented as of the form $Q(L(x,y))$ with $Q$ is an one-variable polynomial and $L$ is a linear form in $x$ and $y$.
\end{theorem}
It also follows from Theorem \ref{vuvu} that  if $|A|=p^\alpha$, then
\[\max\left\lbrace |A+A|, |P(A)|\right\rbrace\gg |A|^{1+\Delta(\alpha)}, \]
where $\Delta(\alpha)=\min(1/2-1/4\alpha, (1/\alpha-1)/3)$.

Recently, Hegyv\'{a}ri and Hennecart \cite{heg} obtained analogous results of these problems by using a generalization of Solymosi's approach in \cite{solymosi}. In particular, they proved that for some certain families of two-variable functions $f(x,y)$ and $g(x,y)$, if $|A|=|B|=p^\alpha$, then $\max\left\lbrace |f(A,B)|, |g(A,B)|\right\rbrace\gg |A|^{1+\Delta(\alpha)}$, for some $\Delta(\alpha)>0$. Before giving their first result, we need the following definition on the multiplicity of a function defined over a subgroup over finite fields.

Let $G$ be a subgroup in $\mathbb{F}_p^*$, and $g\colon G\to \mathbb{F}_p$ an arbitrary function, we define
\[\mu(g)=\max_{t}\left\vert \left\lbrace x\in G\colon g(x)=t\right\rbrace \right\vert.\]
\begin{theorem}[\textbf{Hegyv\'{a}ri and Hennecart}, \cite{heg}]\label{thm1}
Let $G$ be a subgroup of $\mathbb{F}_p^*$, and $f(x,y)=g(x)(h(x)+y)$ be defined on $G\times \mathbb{F}_p^*$, where $g,h\colon G\to \mathbb{F}_p^*$ are arbitrary functions. Put $m=\mu(g\cdot h)$. For any sets $A\subset G$ and $B,C\subset \mathbb{F}_p^*$, we have
\[\left\vert f(A,B)\right\vert \left\vert B\cdot C\right\vert\gg \min\left\lbrace\frac{|A||B|^2|C|}{pm^2}, \frac{p|B|}{m}\right\rbrace.\] 
\end{theorem}
In particular, if $f(x,y)=x(1+y)$, then, as a consequence of Theorem \ref{thm1}, we obtain the following corollary which also studied by Garaev and Shen in \cite{gara}.
\begin{corollary}
For any set $A\subseteq \mathbb{F}_p\setminus \left\lbrace 0, -1\right\rbrace$, we have 
\[|A\cdot (A+1)|\gg \min\left\lbrace\sqrt{p|A|}, |A|^2/\sqrt{p}\right\rbrace.\]
\end{corollary}
The next result is the additive version of Theorem \ref{thm1}.

\begin{theorem}[\textbf{Hegyv\'{a}ri and Hennecart}, \cite{heg}]\label{2.2}
Let $G$ be a subgroup of $\mathbb{F}_p^*$, and $f(x,y)=g(x)(h(x)+y)$ be defined on $G\times \mathbb{F}_q^*$ where $g$ and $h$ are arbitrary functions from $G$ into $\mathbb{F}_p^*$. Put $m=\mu(g)$. For any $A\subset G$, $B,C \subset \mathbb{F}_p^*$, we have
\[|f(A,B)||B+C|\gg \min\left\lbrace  \frac{|A||B|^2|C|}{pm^2}, \frac{p|B|}{m}\right\rbrace\]
\end{theorem}
Note that by letting $C=A$, this implies that 
\[\max\left\lbrace|f(A,B)|,|A+B|\right\rbrace\gg |A|^{1+\Delta(\alpha)}, ~|A|=|B|=p^\alpha,\]
where $\Delta(\alpha)=\min\left\lbrace 1-1/2\alpha, (1/\alpha-1)/2\right\rbrace$. In the case $g$ and $h$ are polynomials, and $g$ is non constant, Theorem \ref{vuvu}, or its generalization in $\cite{bon}$ would lead to a similar statement with a weaker exponent $\Delta(\alpha)=\min\{1/2-1/4\alpha, 1/3\alpha-1/3\}$. We also note that Theorem $6$ established by Bukh and Tsimerman  \cite{bukh} does not cover such a function like in Theorem \ref{2.2}.

For any function $h\colon \mathbb{F}_q\to \mathbb{F}_q$ and $u\in \mathbb{F}_p$, we define $h_u(x):=h(ux)$. In \cite{heg}, Hegyv\'{a}ri and Hennecart obtained a generalization of Theorem \ref{thm1} as follows.
\begin{theorem}[\textbf{Hegyv\'{a}ri and Hennecart}, \cite{heg}]\label{thm4}
Let $f(x,y)=g(x)h(y)(x^k+y^k)$ where $g,h:G\to \mathbb{F}_p^*$ are functions defined on some subgroup $G$ of $\mathbb{F}_p^*$. We assume that for any fixed $z\in G$, $g(xz)/g(x)$ and $h(xz)/h(x)$ take $O(1)$ different values when $x\in G$ and that $\max_u\mu(g\cdot h_u\cdot id)=O(1)$. Then for any $A,B,C\subset G$, one has
\[\left\vert f(A,B) \right\vert|A\cdot C||B\cdot C|\gg \min\left\lbrace \frac{|A|^2|B|^2|C|}{p}, p|A||B|\right\rbrace.\] 
\end{theorem}
The condition on $g$ and $h$ in the theorem looks unusual. For instance, one can take $g$ and $h$ being monomial functions, or functions of the form $\lambda^{\alpha(x)} x^k$, where $\lambda\in \mathbb{F}_p^*$ has order $O(1)$ and $\alpha(x)$ is an arbitrary function. Note that in some particular cases, we can obtain better results. The following theorem is an example.
\begin{theorem}[\textbf{Hegyv\'{a}ri and Hennecart}, \cite{heg}]\label{chat}
Let $A,B,C$ be subsets in $\mathbb{F}_p^*$, and $f(x,y)=xy(x+y)$ a polynomial in $\mathbb{F}_p[x,y]$. Then we have the following estimate
\[|f(A,B)||B\cdot C|\gg \min\left\lbrace \frac{|A||B|^2|C|}{p}, p|B| \right\rbrace.\]
\end{theorem}
This result is sharp when $|A|=|B| \asymp p^\alpha$ with $2/3\le \alpha<1$ since, for instance, one can take $A=B=C$ being a geometric progression of length $p^\alpha$, it is easy to see that $|A\cdot A|\ll |A|$, and $|f(A,A)|\le p$. This implies that $|f(A,A)||A\cdot A|\ll p|A|$.

There is a series of papers dealing with similar results on the sum-product problem, for example, see \cite{mot, hai, ba, bon, batbuoc1, nam, heg1, heg2, batbuoc2, sau, bay, tam}.

Let $\R$ be a finite valuation ring of order $q^r$. Throughout, $\R$ is assumed to be commutative, and to have an identity. Let us denote the set of units, non-units in $\R$ by $\R^{*}, \R^{0}$, respectively. 

The main purpose of this paper is to extend aforementioned results to finite valuation rings by using methods from spectral graph theory. Our first result is a generalization of Theorem \ref{thm1}.
\begin{theorem}\label{thm2}
Let $\R$ be a finite valuation ring of order $q^r$, $G$ be a subgroup of $\R^*$, and $f(x,y)=g(x)(h(x)+y)$ be defined on $G\times \R^*$, where $g,h\colon G\to \R^*$ are arbitrary functions. Put $m=\mu(g\cdot h)$. For any sets $A\subset G$ and $B,C\subset \R^*$, we have
\[|f(A,B)||B\cdot C|\gg \min\left\lbrace \frac{q^r|B|}{m}, \frac{|A||B|^2|C|}{m^2q^{2r-1}}\right\rbrace.\]
\end{theorem}

In the case, $f(x,y)=x(1+y)$, we obtain the following estimate.
\begin{corollary}
For any set $A\subset \R\setminus \{\R^0, \R^0-1\}$, we have
\[\left\vert A(A+1)\right\vert\gg \min\left\lbrace \sqrt{q^r|A|}, \frac{|A|^2}{\sqrt{q^{2r-1}}}\right\rbrace.\]
\end{corollary}
As in Theorem \ref{2.2}, we obtain the additive version of Theorem \ref{thm2} as follows.
\begin{theorem}\label{thm3}
Let $\R$ be a finite valuation ring of order $q^r$, $G$ be a subgroup of $\R^*$, and $f(x,y)=g(x)(h(x)+y)$ be defined on $G\times \R^*$ where $g$ and $h$ are arbitrary functions from $G$ into $\R^*$. Put $m=\mu(g)$. For any $A\subset G$, $B,C \subset \R^*$, we have
\[|f(A,B)||B+ C|\gg \min\left\lbrace \frac{q^r|B|}{m}, \frac{|A||B|^2|C|}{m^2q^{2r-1}}\right\rbrace.\]
\end{theorem}
Combining Theorem \ref{thm2} and Theorem \ref{thm3}, we obtain the following corollary.
\begin{corollary}
Let $f(x,y)=g(x)(x+y)$ such that $\mu(g)=O(1)$, and $A\subset \R^*$. Then 
\[|f(A,A)|\times \min\left\lbrace|A\cdot A|, |A+A|\right\rbrace\gg \min\left\lbrace q^r|A|, \frac{|A|^4}{q^{2r-1}}\right\rbrace.\]
\end{corollary}
Finally, we will derive  generalizations of Theorem \ref{thm4} and Theorem \ref{chat}.
\begin{theorem}\label{thm11} Let $\R$ be a finite valuation ring of order $q^r$, and $f(x,y)=g(x)h(y)(x+y)$ where $g,h:G\to \R^*$ are functions defined on some subgroup $G$ of $\R^*$. We assume that for any fixed $z\in G$, $g(xz)/g(x)$ and $h(xz)/h(x)$ take $O(1)$ different values when $x\in G$ and that $\max_u\mu(g\cdot h_u\cdot id)=O(1)$. Then for any $A,B,C\subset G$, one has
\[\left\vert f(A,B) \right\vert|A\cdot C||B\cdot C|\gg \min\left\lbrace q^r|A||B|, \frac{|A|^2|B|^2|C|}{q^{2r-1}}\right\rbrace.\] 
\end{theorem}
Similarly, we can improve Theorem \ref{thm11} for some special cases of $f(x,y)$. The following theorem is an example, which is  an extension of Theorem \ref{chat}.
\begin{theorem}\label{thm12}
Let $\R$ be a finite valuation ring of order $q^r$, and $A, B, C$ be subsets in $\R^*$, $f(x,y)=xy(g(x)+y)$, where $g$ is a function from $\R^*$ into $\R^*$, and $\mu(g^2\cdot id)=O(1)$. Then we have 
\[|f(A,B)||B\cdot C|\gg \min\left\lbrace q^r|B|, \frac{|A||B|^2|C|}{q^{2r-1}}\right\rbrace.\]
\end{theorem}

Note that we also can obtain similar results over $\mathbb{Z}/m\mathbb{Z}$ by using Lemma 4.1 in \cite{vinh-norm2} instead of Lemma \ref{sp-graph}.

\section{Preliminaries}
We say that a ring $\R$ is \textit{local} if $\R$ has a unique maximal ideal that contains every proper ideal of $\R$. $\R$ is \textit{principal} if every ideal in $\R$ is principal. The following is the definition of finite valuation rings.
\begin{definition}
Finite valuation rings are finite rings that are local and principal.
\end{definition}
Throughout, rings are assumed to be commutative, and to have an identity. Let $\R$ be a finite valuation ring, then $\R$ has a unique maximal ideal that contains every proper ideals of $\R$. This implies that there exists a non-unit $z$ called \textit{uniformizer} in $\R$ such that the maximal ideal is generated by $z$. Throughout this paper, we denote the maximal ideal of $\R$ by $(z)$. Moreover, we also note that the uniformizer $z$ is defined up to a unit of $\R$. 

There are two structural parameters associated to $\R$ as follows: the cardinality of the residue field $F=\R/(z)$, and the nilpotency degree of $z$, where the nilpotency degree of $z$ is the smallest integer $r$ such that $z^r=0$.  Let us denote the cardinality of $F$ by $q$. 
In this note, $q$ is assumed to be odd, then $2$ is a unit in $\R$.

If $\R$ is a finite valuation ring, and $r$ is the nilpotency degree of $z$, then we have a natural valuation
\[\nu\colon \R\to \{0,1,\ldots,r\}\]
defined as follows: $\nu(0)=r$, for $x\ne 0$, $\nu(x)=k$ if $x\in (z^k)\setminus (z^{k+1})$. We also note that $\nu(x)=k$ if and only if $x=uz^k$ for some unit $u$ in $\R$. Each abelian group $(z^k)/(z^{k+1})$ is a one-dimensional linear space over the residue field $F=\R/(z)$, thus its size is $q$. This implies that $|(z^k)|=q^{r-k}, ~k=0,1,\ldots,r$. In particular, $|(z)|=q^{r-1}, |\R|=q^r$ and $|\R^{*}|=|\R|-|(z)|=q^r-q^{r-1}$, (for more details about valuation rings, see \cite{ati}, \cite{bi}, \cite{fulton}, and \cite{bogan}). The following are some examples of finite valuation rings:
\begin{enumerate}
\item Finite fields $\mathbb{F}_q$, $q=p^n$ for some $n>0$.
\item Finite rings $\mathbb{Z}/p^r\mathbb{Z}$, where $p$ is a prime.
\item $\mathcal{O}/(p^r)$ where $\mathcal{O}$ is the ring of integers in a number field and $p\in \mathcal{O}$ is a prime.
\item $\mathbb{F}_q[x]/(f^r)$, where $f\in \mathbb{F}_q[x]$ is an irreducible polynomial.
\end{enumerate}
\section{Properties of pseudo-random graphs}

For a graph $G$ of order $n$, let $\lambda_1 \geq \lambda_2 \geq \ldots \geq \lambda_n$ be
the eigenvalues of its adjacency matrix. The quantity $\lambda (G) = \max
\{\lambda_2, - \lambda_n \}$ is called the second eigenvalue of $G$. A graph $G
= (V, E)$ is called an $(n, d, \lambda)$-graph if it is $d$-regular, has $n$
vertices, and the second eigenvalue of $G$ is at most $\lambda$. Since $G$ is a $d$-regular graph, $d$ is an eigenvalue of its adjacency matrix with the all-one eigenvector $\tmmathbf{1}$. If the graph $G$ is connected, the eigenvalue $d$ has multiplicity one. Furthermore, if $G$ is not bipartite, for any other eigenvalue $\theta$ of $G$, we have $|\theta| < d$. Let $\tmmathbf{v}_{\theta}$ denote the corresponding eigenvector of $\theta$. We will make use of the trick that $\tmmathbf{v}_\theta \in \tmmathbf{1}^{\bot}$, so $J\tmmathbf{v}_\theta = 0$ where $J$ is the all-one matrix of size $n \times n$ (see \cite{bh} for more background on spectral graph theory).

It is well-known (see \cite[Chapter 9]{as} for more details) that if $\lambda$ is much smaller than the
degree $d$, then $G$ has certain random-like properties.  For two (not necessarily) disjoint subsets of vertices $U,
W \subset V$, let $e (U, W)$ be the number of ordered pairs $(u, w)$ such that
$u \in U$, $w \in W$, and $(u, w)$ is an edge of $G$. We recall the following 
well-known fact (see, for example, \cite{as}).

\begin{lemma}(\cite[Corollary 9.2.5]{as})\label{edge}
  Let $G = (V, E)$ be an $(n, d, \lambda)$-graph. For any two sets $B, C
  \subset V$, we have
  \[ \left| e (B, C) - \frac{d|B | |C|}{n} \right| \leq \lambda \sqrt{|B| |C|}. \]
\end{lemma}

\subsection{Sum-product graphs over finite valuation rings}
The sum-product (undirected) graph $\mathcal{SP}_\R$ is defined as follows. The vertex set of the
sum-product graph $\mathcal{S P}_\R$ is the set $V (\mathcal{S P}_\R)
=\R \times \R$. Two vertices $U = (a, b)$ and
$V = (c, d) \in V (\mathcal{SP}_\R)$ are connected by an edge, $(U, V)
\in E (\mathcal{SP}_\R)$, if and only if $a + c= bd$. Our construction is similar to that of Solymosi in \cite{solymosi}.

\begin{lemma}\label{sp-graph} Let $\R$ be a finite valuation ring. The sum-product graph, $\mathcal{SP}_\R$, is a \[\left(q^{2r}, q^{r},\sqrt{2rq^{2r-1}}\right)-\mbox{graph.}\]
\end{lemma}

\begin{proof}
  It is easy to see that $\mathcal{SP}_\R$ is a regular graph of order
  $q^{2r}$ and valency $q^{r}$. We now compute the eigenvalues of this
  multigraph (there are few loops). For any two vertices $(a, b), (c, d)  \in \R\times\R$, we count the number of solutions of the following system
  \begin{equation} \label{esps}a + u  = bv,\,\, c + u =
    dv,\,\, \; (u,v) \in \R\times \R. \end{equation}
For each solution $v$ of \begin{equation}\label{esps1}(b-d)v= a - c,\end{equation} there exists a
  unique $u$ satisfying the system (\ref{esps}). Therefore, we only need to count the number of solutions of (\ref{esps1}). Suppose that $\nu(b-d)=\alpha$. If $\nu(a-c)< \alpha$, then Eq. (\ref{esps1}) has no solution. Thus we assume that $\nu(a-c)\ge \alpha$. It follows from the definition of the function $\nu$ that there exist $u_1, u_2$ in $\R^*$ such that $a-c=u_1z^{\nu(a-c)}, ~b-d=u_2z^{\nu(b-d)}$.  Let $\mu = u_1z^{\nu(a-c)-\alpha}$ and $x= u_2z^{\nu(b-d)-\alpha}$. The number of solutions of (\ref{esps1}) equals the number of solutions $v\in \R$ satisfying 
  \begin{equation}\label{eq10}x\cdot v-\mu\in (z^{r-\alpha}).\end{equation}
Since $\nu(b-d)=\alpha$, we have $x\in \R^*$, and the equation 
\[xv-\mu=t\]
has a unique solution for each $t\in (z^{r-\alpha})$. Since $|(z^{r-\alpha})|=q^\alpha$, the number solutions of (\ref{eq10}) is $q^{\alpha}$  if $\nu(a-c)\ge \alpha$. 
  
Therefore, for any two vertices $U =
  (a, b)$ and $V = (c, d) \in V (\mathcal{SP}_\R)$,  $U$ and $V$ have $q^{\alpha}$ common neighbors if $\nu(b-d)=\alpha$ and  $\nu(a-c)\ge \alpha$ and no common neighbor if $\nu(b-d)=\alpha$ and  $\nu(c- a)< \alpha$. Let $A$ be the adjacency matrix of $\mathcal{SP}_\R$. For any two vertices $U, V$ then $(A^2)_{U,V}$ is the number of common vertices of $U$ and $V$. It follows that
  \begin{equation}\label{e20-sesb}
 A^2 = J + (q^{r}-1)I - \sum_{\alpha=0}^{r}E_\alpha + \sum_{\alpha=1}^{r-1}(q^{\alpha}-1)F_\alpha,
 \end{equation}
  where:
\begin{itemize}
\item $J$ is the all-one matrix and $I$ is the identity matrix.
\item $E_\alpha$ is the adjacency matrix of the graph $B_{E,\alpha}$, where for any two vertices $U =
  (a, b)$ and $V = (c, d) \in V (\mathcal{SP}_\R)$, $(U,V)$ is an edge of $B_{E,\alpha}$ if and only if  $\nu(b-d)=\alpha$ and  $\nu(a-c)<\alpha$
\item $F_\alpha$ is the adjacency matrix of the graph $B_{F,\alpha}$, where for any two vertices $U =
  (a, b)$ and $V = (c, d) \in V (\mathcal{SP}_\R)$, $(U,V)$ is an edge of $B_{F,\alpha}$ if and only if  $\nu(b-d)=\alpha$ and  $\nu(a-c)\ge \alpha$
\end{itemize}

For any $\alpha> 0$, we have $|(z^\alpha)|=q^{r-\alpha}$, thus $B_{E,\alpha}$ is a regular graph of valency less than $q^{2r-\alpha}$ and $B_{F,\alpha}$ is a regular graph of valency less than $q^{2(r-\alpha)}$. Since eigenvalues of a regular graph are bounded by its valency, all eigenvalues of $E_\alpha$ are at most  $q^{2r-\alpha}$ and all eigenvalues of $F_\alpha$ are at most $q^{2(r-\alpha)}$. Note that $E_0$ is a zero matrix. 
  
  Since $\mathcal{SP}_\R$ is a $q^{r}$-regular graph, $q^{r}$ is an eigenvalue of $A$ with the all-one eigenvector $\tmmathbf{1}$. The graph $\mathcal{SP}_\R$ is connected therefore the eigenvalue $q^{r}$ has multiplicity one. Note that for two adjacent vertices $U=(2z^{2\alpha+1}, z^\alpha)$ and $V=(-z^{2\alpha+1}, z^{\alpha+1})$, they have many common neighbors. This implies that the graph $\mathcal{SP}_\R$ contains (many) triangles, it is not bipartite. In the case $|(z)|=1$, then $U=V$, and $\R$ is a finite field, we can also check that it contains many triangles. Hence, for any other eigenvalue $\theta$, $|\theta| < q^{r}$. Let $\tmmathbf{v}_\theta$ denote the corresponding eigenvector of $\theta$. Note that $\tmmathbf{v}_\theta \in \tmmathbf{1}^{\bot}$, so $J\tmmathbf{v}_\theta$ = 0. It follows from (\ref{e20-sesb}) that
\[(\theta^2 - q^{r}+ 1)\tmmathbf{v}_\theta = \left(\sum_{\alpha=1}^{r}E_\alpha - \sum_{\alpha=1}^{r-1}(q^{\alpha}-1)F_\alpha\right)\tmmathbf{v}_\theta.\]
Hence, $\tmmathbf{v}_\theta$ is also an eigenvalue of  \[\sum_{\alpha=1}^{r}E_\alpha - \sum_{\alpha=1}^{r-1}(q^{\alpha}-1)F_\alpha\] Since absolute value of eigenvalues of sum of matrices are bounded by sum of largest absolute values of eigenvalues of summands. We have
\begin{eqnarray} \theta^2 &\leq& q^{r}-1 + \sum_{\alpha=1}^{r}q^{2r-\alpha} + \sum_{\alpha=1}^{r-1}(q^{\alpha}-1)q^{2(r-\alpha)}\nonumber\\
& < & 2rq^{2r-1}.\nonumber
\end{eqnarray}
The lemma follows.
\end{proof}

\section{Proofs of Theorem \ref{thm2} and Theorem \ref{thm3}}
\begin{proof}[Proof of Theorem \ref{thm2}]
First we set
\[S=\left\lbrace \left( zh(x), zg(x)^{-1}\right)\colon (x,z)\in A\times C  \right\rbrace, T=\left\lbrace \left( yz, g(x)(h(x)+y)\right)\colon (x,y,z)\in A\times B\times C  \right\rbrace\]
This implies that 
\[|S|\le |A||C|, |T|\le \min\left\lbrace |A||B||C|, |f(A,B)||B\cdot C|\right\rbrace.\]
Given a quadruple $(u,v,w,t)\in (\R^*)^4$, we now count the number of solutions $(x,y,z)$ to the following system
\[g(x)(h(x)+y)=u,~ yz=v, ~zg(x)^{-1}=w, ~zh(x)=t.\]
This implies that 
\[ g(x)h(x)=\frac{t}{w}=\frac{ut}{v+t}.\]
Since $\mu(g\cdot h)=m$, there are at most $m$ different values of $x$ satisfying the equality $g(x)h(x)=t/w$, and $y, z$ are determined uniquely in terms of $x$ by the second and the fourth equations. Therefore, the number of edges between $S$ and $T$ in the sum-product graph $\mathcal{SP}_\R$ is at least $|A||B||C|/m$. On the other hand, it follows from Lemma \ref{edge} and Lemma \ref{sp-graph} that 
\[\frac{|A||B||C|}{m}\le e(S,T)\le \frac{|S||T|}{q^r}+\sqrt{2r}q^{(2r-1)/2}\sqrt{|S||T|}.\]
Solving this inequality gives us 
\[|S||T|\gg \min\left\lbrace \frac{q^{r}|A||B||C|}{m}, \frac{(|A||B||C|)^2}{m^2q^{2r-1}}\right\rbrace.\]
Thus, we obtain
\[|f(A,B)||B\cdot C|\gg \min\left\lbrace \frac{q^r|B|}{m}, \frac{|A||B|^2|C|}{m^2q^{2r-1}}\right\rbrace,\]
which concludes the proof of theorem.\end{proof}
\begin{proof}[Proof of Theorem \ref{thm3}]
The proof of Theorem \ref{thm3} is as similar as  the proof of Theorem \ref{thm2} by setting
\[S=\left\lbrace (y+z, g(x)(h(x)+y))\colon (x,y,z)\in A\times B\times C\right\rbrace,\]
\[T=\left\lbrace (h(x)-z, g(x)^{-1})\colon (x,y,z)\in A\times B\times C\right\rbrace.\]
\end{proof}

\section{Proofs of Theorem \ref{thm11} and Theorem \ref{thm12}}
\begin{proof}[Proof of Theorem \ref{thm11}]
Let
\[S=\left\lbrace \left(yz, \frac{g(x)h(y)(x+y)}{h(yz)}\right)\colon (x,y,z)\in A\times B\times C\right\rbrace,\]
\[T=\left\lbrace \left(xz, \frac{zg(xz)h(yz)g(x)^{-1}h(y)^{-1}}{g(xz)}\right)\colon (x,y,z)\in A\times B\times C\right\rbrace.\]
Then $S$ and $T$ are two sets of vertices in the sum-product graph $\mathcal{SP}_\R$, and $|S|\ll|f(A,B)||B\cdot C|$, $|T|\ll |C||A\cdot C|$.
Given a quadruple $(u,v,w,t)$ in $(\R^*)^4$, we now count the number of solutions $(x,y,z)$ to the following system
\[\frac{g(x)h(y)(x+y)}{h(yz)}=u, ~yz=v, ~\frac{zg(xz)h(yz)g(x)^{-1}h(y)^{-1}}{g(xz)}=t, ~zx=w.\]
This implies that \begin{equation}\label{eq100}xg(x)h(vx/w)=\frac{uw}{w+v}h(v).\end{equation}
Since $\max_u \mu(g\cdot h_u\cdot id)=O(1)$, there are at most $O(1)$ values of $x$ satisfying the equation (\ref{eq100}), and $y,z$ are determined uniquely in terms of $x$ by the second and the fourth equations. Thus, the number of edges between $S$ and $T$ in $\mathcal{SP}_\R$ is at least $\gg |A||B||C|$. The rest of the proof is the same as the proof of Theorem \ref{thm2}.
\end{proof}

\begin{proof}[Proof of Theorem \ref{thm12}]
First we set
\[S=\left\lbrace \left(yz, \frac{xy(g(x)+y)}{yz}\right)\colon (x,y,z)\in A\times B\times C\right\rbrace, ~T=\left\lbrace \left(zg(x), \frac{z^{2}}{x} \right)\colon (x,z)\in A\times C\right\rbrace.\]
Then $S$ and $T$ are two sets of vertices in the sum-product graph $SP_\R$, and  $|S|\le |f(A,B)||B\cdot C|$, $|T|\le |A||C|$. It follows from Lemma \ref{edge} and Lemma \ref{sp-graph} that 
\begin{equation}\label{eq90}
e(S,T)\le \frac{|S||T|}{q^r}+\sqrt{2r}q^{(2r-1)/2}\sqrt{|S||T|}.\end{equation}
On the other hand, given a quadruple $(u,v,w,t)$ in $(\R^*)^4$, we now count the number of solutions $(x,y,z)$ to the following system
\[\frac{xy(g(x)+y)}{yz}=u, ~yz=v, ~\frac{z^2}{x}=t, ~zg(x)=w.\]
This implies that $g(x)^2x=w^2/t$. Since $\mu(g^2\cdot id)=O(1)$, there are at most $O(1)$ values of $x$ satisfying the equality $g(x)^2x=w^2/t$, and $y, z$ are determined uniquely in terms of $x$ by the second and the fourth equations. Therefore, we have
\begin{equation}\label{eq91}e(S,T)\gg |A||B||C|.\end{equation}
Putting (\ref{eq90}) and (\ref{eq91}) together, we get
\[|A||B||C|\ll \frac{|S||T|}{q^r}+\sqrt{2r}q^{(2r-1)/2}\sqrt{|S||T|}.\]
This implies that 
\[|S||T|\gg \min\left\lbrace  q^{r}|A||B||C|, \frac{(|A||B||C|)^2}{q^{2r-1}}\right\rbrace.\]
Therefore,
\[|f(A,B)||B\cdot C|\gg \min\left\lbrace q^r|B|, \frac{|A||B|^2|C|}{q^{2r-1}}\right\rbrace,\]
and the theorem follows.
\end{proof}
\section*{Acknowledgements.}
\thispagestyle{empty}
The authors would like to thank two anonymous referees for valuable comments and suggestions which improved the presentation of this paper considerably.


\begin{thebibliography}{99}
\bibitem{as}
N. Alon and J. H. Spencer, \textit{The probabilistic method}, 2nd ed., Willey-Interscience, 2000.
\bibitem{ati}
M.F. Atiyah, I.G. Macdonald, \textit{Introduction to commutative algebra} (Vol. 2). Reading: Addison-Wesley.(1969)

\bibitem{bourgain-katz-tao}
J. Bourgain, N. Katz, T. Tao, \textit{A sum-product estimate in finite fields, and applications}, Geom. Funct. Anal. \textbf{14} (2004), 27--57.
\bibitem{mot}
J. Bourgain, \textit{More on the sum-product phenomenon in prime fields and its applications}, International Journal of Number Theory, \textbf{1}(01) (2005), 1--32.
\bibitem{hai}
J. Bourgain, M. Z. Garaev, \textit{On a variant of sum-product estimates and explicit exponential sum
bounds in prime fields}, Math. Proc. Cambridge Philos. Soc. \textbf{146} (2009), no. 1, 1--21.

\bibitem{bukh}
B. Bukh,  J. Tsimerman, \textit{Sum–product estimates for rational functions}, Proceedings of the London Mathematical Society, (2011), \textbf{104} (1), 1--26.
\bibitem{bh} A. Brouwer and W. Haemers, \textit{Spectra of Graphs}, Springer, New York, etc., 2012. 
\bibitem{bi}
G. Bini, F. Flamini, \textit{Finite commutative rings and their applications}, Kluwer International Series in Engineering and Computer Science 680, Kluwer Academic Publishers 2002.
\bibitem{cil}
J.\ Cilleruelo, \textit{Combinatorial problems in finite fields and Sidon sets}, Combinatorica {\bf 32}(5) (2012), 497--511.
\bibitem{fulton}
W. Fulton, \textit{Algebraic curves: An introduction to algebraic geometry},  Notes written
with the collaboration of Richard Weiss. Reprint of 1969 original. Advanced Book
Classics. Addison-Wesley Publishing Company, Advanced Book Program, Redwood
City, CA, (1989).
\bibitem{garaev}
M.\ Z.\ Garaev, \textit{The sum-product estimate for large subsets of prime fields},
Proc.\ Amer.\ Math.\ Soc., \textbf{136}(2008), 2735--2739.
\bibitem{gara}
M. Garaev, C.-Y. Shen, \textit{On the size of the set $A(A+1)$}, Math. Z. \textbf{263}(2009), no. 94.
\bibitem{ba}
A. A. Glibichuk,  S. V. Konyagin, \textit{ Additive properties of product sets in prime fields order, Additive
combinatorics}, 279–286, CRM Proc. Lecture Notes, 43, Amer. Math. Soc., Providence, RI,
(2007).
\bibitem{bon}
D. Hart, L. Li, C-Y. Shen, \textit{Fourier analysis and expanding phenomena in finite fields}, Proceedings of the American Mathematical Society, \textbf{141}(2)(2013), 461--473.
\bibitem{his}
D.\ Hart, A.\ Iosevich, J.\ Solymosi, \textit{Sum-product estimates in finite fields via Kloosterman sums}, Int.\ Math.\ Res.\ Not.\  no. 5, (2007) Art.\ ID rnm007.   
\bibitem{batbuoc1}
H. A. Helfgott, M. Rudnev, \textit{An explicit incidence theorem in $\mathbb{F}_p$}, Mathematika, \textbf{57} (2011) 135--145.
\bibitem{nam}
N. Hegyv\'{a}ri, F. Hennecart, \textit{Explicit construction of extractors and expanders}, Acta Arith. \textbf{140}(2009), 233--249.
\bibitem{heg}
N. Hegyv\'{a}ri, F. Hennecart, \textit{Conditional expanding bounds for two-variable functions
over prime fields}, European J. Combin., \textbf{34}(2013), 1365--1382.
\bibitem{heg1}
N. Hegyv\'{a}ri,  \textit{Some remarks on multilinear exponential sums with an application}, Journal of Number Theory, \textbf{132}(1) (2012), 94--102.
\bibitem{heg2}
N. Hegyv\'{a}ri, F. Hennecart,  \textit{A structure result for bricks in Heisenberg groups}, Journal of Number Theory, \textbf{133}(9) (2013), 2999--3006.
\bibitem{batbuoc2}
T. Jones, \textit{An improved incidence bound for fields of prime order}, European Journal of Combinatorics, Volume 52, Part A, February 2016, Pages 136--145
\bibitem{sau}
N. H. Katz, C-Y. Shen, \textit{A slight improvement to Garaev's sum product estimate}, Proc. Amer.
Math. Soc. \textbf{136} (2008), 2499--2504.
\bibitem{bogan}
B. Nica,  \textit{Unimodular graphs and Eisenstein sums}, arXiv: 1505.05034 (2015).
\bibitem{bay}
L. Li, \textit{Slightly improved sum-product estimates in fields of prime order}, Acta Arith. \textbf{147} (2011),no. 2, 153--160.
\bibitem{solymosi}J. Solymosi, \textit{Incidences and the Spectra of Graphs}, Building Bridges between Mathematics and Computer Science. Vol. \textbf{19.} Ed. Martin Groetschel and Gyula Katona. Series: Bolyai Society Mathematical Studies. Springer (2008), 499 -- 513.
\bibitem{tam}
T. Tao, \textit{Expanding polynomials over finite fields of large characteristic, and a regularity lemma for definable sets}, Contribution to Discrete Mathematics in Volume 10, Number 1, Pages 22--98.
\bibitem{vinh}
L.A.Vinh, \textit{A Szemer\'{e}di--Trotter type theorem and sum-product estimate over finite fields}, European J. Combin.,  \textbf{32}(2011), no. 8, 1177--1181.
\bibitem{vinh-norm2}
L.A.Vinh, Product graphs, Sum-product graphs and sum-product estimate over finite rings, \textit{Forum Mathematicum}, Volume 27, Issue 3 (2015), 1639--1655.
\bibitem{vu}
H. V. Vu, \textit{Sum-product estimates via directed expanders}, Mathematical research letters, \textbf{15}(2) (2008), 375--388.
\end{thebibliography}
\end{document}